\documentclass[amsfonts,11pt]{amsart}
\RequirePackage{amsmath,amssymb,amsthm, amscd, comment,mathtools}
\usepackage[margin=1.2in]{geometry}
\usepackage{amsmath}
\usepackage{amsthm}
\usepackage{amssymb}
\usepackage{amscd}
\usepackage{mathrsfs}
\usepackage{graphicx}
\usepackage{rotating, graphicx}
\usepackage[all,cmtip]{xy}
\usepackage{enumerate}
\usepackage{tabularx}
\usepackage{tikz-cd}
\usepackage[thinlines]{easytable}
\usepackage{multirow}
\usepackage{mathtools}
\usepackage{mathrsfs}
\usepackage{changepage} 
\usepackage{paralist}

\usepackage{kantlipsum}
\usepackage{mathabx}
\usepackage{enumitem}
\usepackage{times}
\usepackage{subfig}
\usepackage{color}
\usepackage[pagebackref,breaklinks,colorlinks,linkcolor=red,anchorcolor=red,citecolor=blue]{hyperref}

\calclayout

\newtheorem{proposition}[subsection]{Proposition}
\newtheorem{corollary}[subsection]{Corollary}
\newtheorem{theorem}[subsection]{Theorem}
\newtheorem{lemma}[subsection]{Lemma}
\theoremstyle{definition}
\newtheorem{definition}[subsection]{Definition}
\theoremstyle{remark}

\newtheorem{remark}[subsection]{Remark}

\numberwithin{equation}{subsection}

\DeclareMathAlphabet{\mathbbold}{U}{bbold}{m}{n}

\title{On some finiteness results in real \'etale cohomology}

\author{Fangzhou Jin}
\address{School of Mathematical Sciences\\
Tongji University\\
Siping Road 1239\\
200092 Shanghai\\
China}
\address{Fakult\"at f\"ur Mathematik,
Universit\"at Duisburg-Essen,
Thea-Leymann-Strasse 9,
45127 Essen,
Germany}
\email{\href{mailto:fangzhoujin@tongji.edu.cn}{fangzhoujin@tongji.edu.cn}}
\urladdr{\url{https://fangzhoujin.github.io/}}

\date{\number\day-\number\month-\number\year}

\begin{document}

\maketitle

\begin{abstract}

We show that for quasi-compact quasi-separated schemes of finite dimension, the constructibility condition in real \'etale cohomology agrees with a notion of constructibility arising naturally from topology. As application we prove that the derived direct image functor preserves constructibility under some assumptions, and compute the Grothendieck group of the constructible rational stable motivic homotopy category for reasonable schemes.
We prove the generic base change property for constructible real \'etale sheaves, and deduce the same property for rational motivic spectra and $b$-sheaves.

\end{abstract}

\setcounter{tocdepth}{1}
\tableofcontents

\noindent

\section{Introduction}

\subsection{}
The study of \emph{real \'etale cohomology} originates from attempts to understand the link between \'etale cohomology with $2$-torsion coefficients and real algebraic geometry. To any scheme $X$, one can associate a \emph{real scheme} $X_r$ which is a topological space obtained by gluing the real spectra of rings (\cite[0.4.2]{Sch3}). The points of $X_r$ are pairs $(x,\sigma)$, where $x$ is a point of $X$ and $\sigma$ is an ordering on the residue field of $x$. The \emph{real \'etale site} of $X$ is the category of \'etale $X$-schemes endowed with the Grothendieck topology where covering families are the ones that induce surjections on real schemes (\cite{CRC}, \cite[Definition 1.2.1]{Sch3}). A fundamental theorem of Coste-Roy-Coste and Scheiderer (\cite[Theorem 1.3]{Sch3}) states that there is a canonical equivalence of sites between the real \'etale site on $X$ and the site defined by the topological space $X_r$, and in particular the real \'etale site is spatial. This result makes it much convenient when working on real \'etale cohomology, as one can not only find inspirations from the powerful machinery of \'etale cohomology (see \cite{SGA4.5}), but also one benefit from the fact that it is possible to work with an actual topological space instead of an abstract Grothendieck topology.

\subsection{}
The following deep result shows an intrinsic relation between the real \'etale sheaves and motivic homotopy theory: for any noetherian scheme of finite dimension $X$, there is a canonical equivalence
\begin{align}
\label{eq:SHdec}
\mathbf{SH}(X,\mathbb{Q})
\simeq
\mathbf{DM}(X,\mathbb{Q})
\times
D(X_r,\mathbb{Q}).
\end{align}
where $\mathbf{SH}(X,\mathbb{Q})$ is the rational stable motivic homotopy category over $X$, $\mathbf{DM}(X,\mathbb{Q})$ is the category of Beilinson motives over $X$ (\cite[Definition 14.2.1]{CD2}), and $D(X_r,\mathbb{Q})$ is the derived category of sheaves over $X_r$.\footnote{More generally, it is proved in \cite{Bac} that real \'etale motivic stable homotopy category is equivalent to $\rho$-inverted motivic stable homotopy category. As the results of this paper only concern the rational case, we only state this version.}

Indeed, by a result of Morel, the ``switching factors'' automorphism of $\mathbb{P}^1\wedge\mathbb{P}^1$ induces a decomposition of $\mathbf{SH}(X,\mathbb{Q})$ into the $+$-part and the $-$-part; the $+$-part is identified with $\mathbf{DM}(X,\mathbb{Q})$ by \cite[Theorem 16.2.13]{CD2}, and the $-$-part is identified with the derived category $D(X_r,\mathbb{Q})$ by \cite[Theorem 35]{Bac} combined with the equivalence between $\mathbf{SH}$ and $D^{\mathbb{A}^1}$ with rational coefficients.

Furthermore, it is proved in \cite{Bac} that the association $X\mapsto D(X_r,\Lambda)$ for a ring $\Lambda$ satisfies the axioms of a \emph{motivic $\infty$-category} (\cite{Kha}),\footnote{One may readily replace this notion with the one of motivic triangulated categories in \cite[Definition 2.4.45]{CD2}.} and therefore one can define and study the associated six functors.\footnote{The results in \cite{Bac} are stated for noetherian schemes of finite dimension, which can be extended to quasi-compact quasi-separated schemes using continuity arguments, at least in the case of the derived category, see \cite[Appendix A]{DFKJ} and \cite[Lemma B.10]{ES}.} Another treatment of real \'etale motivic homotopy theory from the viewpoint of $C_2$-equivariant homotopy theory can be found in \cite{ES}.

\subsection{}
\label{num:motcons}
The goal of this paper is to investigate some finiteness conditions on $D(X_r,\Lambda)$ using these structures. In Section~\ref{sec:cons} we study constructibility conditions. First there is a natural constructibility condition on motivic categories introduced by Ayoub (\cite[D\'efinition 2.2.3]{Ayo}): if $\mathcal{T}$ is a motivic $\infty$-category and $X$ is a scheme, the subcategory of \textbf{constructible objects} over $X$, denoted as $\mathcal{T}_c(X)$, is the thick subcategory generated by elements of the form $Rf_!f^!\mathbbold{1}_X(d)$, where $f:Y\to X$ is a smooth morphism and $(d)$ stands for the Tate twist. This notion can be readily recognized categorically: if $\mathcal{T}$ satisfies some compact generation properties, then the subcategory of constructible objects are exactly the compact objects (see Lemma~\ref{lm:mconscomp}). In particular, the equivalence~\eqref{eq:SHdec} induces an equivalence of the subcategories of constructible objects.

On the other hand, the real scheme associated to a quasi-compact quasi-separated scheme is a \emph{spectral space} (see~\ref{num:spectral}), and there is a natural notion of constructible sheaves on spectral spaces (\cite[Definition A.3]{Sch3}). We promote this notion to the derived category in Definition~\ref{def:conscom} (see also \cite[\S20]{Sch}). The first main result is that under some assumptions these two notions of constructibility agree:
\begin{theorem}[See Theorem~\ref{th:ctf=c}]
Let $X$ be a quasi-compact quasi-separated scheme of finite dimension and let $\Lambda$ be a ring. Then the constructible objects in $D_c(X_r,\Lambda)$ are exactly the complexes $C$ such that there exists a finite stratification of $X_r$ into locally closed constructible subsets such that the restriction of $C$ to each stratum is the constant sheaf associated to a perfect complex of $\Lambda$-modules.
\end{theorem}
The proof reduces to showing that the notions of constructibility and compactness are equivalent, which holds more generally for spectral spaces of finite dimension (see Proposition~\ref{prop:conscomp}): the key result is a theorem of Scheiderer (\cite[Corollary 4.6]{Sch2}) which states that the cohomological dimension of a spectral space is bounded by the Krull dimension. Using this result we are able to give a partial answer to a question raised by Scheiderer (\cite[Remark 17.7.1]{Sch3}), by proving that under some assumptions the derived direct image functor preserves constructible objects, see Corollary~\ref{cor:prescons} and Remark~\ref{rm:schfin}. Another application is the following description of the Grothendieck group of the constructible rational stable motivic homotopy category:
\begin{theorem}[See Corollary~\ref{cor:K0SHQ}]
For $X$ an excellent separated scheme of finite dimension, there is a canonical isomorphism of abelian groups
\begin{align}
K_0(\mathbf{SH}_c(X,\mathbb{Q}))\simeq K_0^{\oplus}(\mathbf{CHM}(X,\mathbb{Q}))\oplus Cons(X_r,\mathbb{Z}).
\end{align}
\end{theorem}
Here $K_0^{\oplus}(\mathbf{CHM}(X,\mathbb{Q}))$ is the direct sum Grothendieck group of the additive category of Chow motives  over $X$ with rational coefficients (see~\ref{num:CHM}), and $Cons(X_r,\mathbb{Z})$ is the group of $\mathbb{Z}$-valued constructible functions on the spectral space $X_r$ (see~\ref{num:consfct}), which is a free abelian group generated by the closed constructible subsets of $X_r$. The proof uses Bondarko's theory of weight structures (see~\ref{num:CHM}), and a general comparison result between the Grothendieck group of the constructible complexes and the group of constructible functions over a spectral space (Proposition~\ref{prop:KS9.7.1}). Note that there are also analogous results for the category $\mathbf{SH}_c(X_r)$, see Remark~\ref{rm:SHret}.

\subsection{}
In Section~\ref{sec:genbc} we prove the generic base change property (see~\ref{num:gbc}) for constructible real \'etale sheaves. The generic base change theorem for \'etale sheaves is due to Deligne (\cite[Th. finitude 1.9]{SGA4.5}), and is generalized to $h$-motives in \cite[2.4.2]{Cis}. Our treatment here uses a very similar strategy, together with some inputs from the topology of real schemes. In Theorem~\ref{th:ret_genbc} we prove the generic base change property for constructible complexes of real \'etale sheaves. From this we deduce the generic base change property for constructible rational motivic spectra (Corollary~\ref{cor:genbcsh}) and for constructible complexes of $b$-sheaves (Corollary~\ref{cor:genbcb}).

\subsection{Conventions and notations}
All smooth morphisms of schemes are assumed separated of finite presentation, and the dimension of a topological space or a scheme stands for the Krull dimension. 

We use the following notations: if $\mathcal{T}$ is a fibered category, $f:X\to Y$ is a morphism and $K\in\mathcal{T}(Y)$, we denote $K_{|X}=f^*K\in\mathcal{T}(X)$. If $f:X\to Y$ is a morphism of schemes, we denote $f^*:D(Y_r,\Lambda)\to D(X_r,\Lambda)$ instead of $f_r^*$, and similarly for the other functors, to simplify the notations.

\subsection{}
\label{num:char0}
For any scheme $X$ we have $X_r=(X\times_\mathbb{Z}\mathbb{Q})_r$, that is, the underlying real scheme only depends on the characteristic $0$ fiber: this is because any field of positive characteristic cannot have an ordering. Therefore in the study of real schemes it is harmless to assume that all schemes have characteristic $0$.

\subsection{}
The following localization sequence is standard in derived categories (see for example \cite[IV Proposition 14.6]{SGA4}): let $X$ be a topological space, let $i:Z\to X$ be the inclusion of a closed subspace with open complement $j:U\to X$ and let $\Lambda$ be a ring, then there is a canonical cofiber sequence
\begin{align}
\label{eq:locseq}
Rj_!j^*\to id\to Ri_*i^*\to Rj_!j^*[1]
\end{align}
of endofunctors of the derived category $D(X,\Lambda)$.

\subsection*{Acknowledgments}
The author would like to thank Denis-Charles Cisinski, Fr\'ed\'eric D\'eglise, Jean Fasel, Adeel Khan, Marc Levine and Heng Xie for helpful discussions. He would also like to thank the referee for carefully reading the paper and for the comments and suggestions which helped improving the quality of the paper. He acknowledges support from the National Key Research and Development Program of China Grant Nr.2021YFA1001400, the National Natural Science Foundation of China Grant Nr.12101455, the Fundamental Research Funds for the Central Universities, and the ERC Project-QUADAG, which has received funding from the European Research Council (ERC) under the European Union’s Horizon 2020 research and innovation programme (grant agreement Nr.832833).

\section{Constructibility}
\label{sec:cons}
In this section we discuss constructibility conditions in $D(X_r,\Lambda)$, and deduce some consequences on Grothendieck groups. Throughout the section, we assume that $\Lambda$ is a commutative ring.
\subsection{}
We first look at the motivic constructibility condition as recalled in~\ref{num:motcons}. The following result is standard:
\begin{lemma}
\label{lm:mconscomp}
Let $X$ be quasi-compact quasi-separated scheme. Then $D(X_r,\Lambda)$ is compactly generated, and the objects of $D_c(X_r,\Lambda)$ are exactly the compact objects of $D(X_r,\Lambda)$.
\end{lemma}
\proof
Since $D(X_r,\Lambda)$ is compactly generated by the Tate twists by \cite[Lemma 15 and Theorem 35]{Bac}, the result follows from \cite[Proposition 1.4.11]{CD2}.
\endproof

\subsection{}
\label{num:spectral}
We now consider a topological constructibility condition. 
Recall that a \textbf{spectral space} is a quasi-compact topological space in which
\begin{enumerate}
\item the intersection of two quasi-compact opens is quasi-compact (i.e. the space is quasi-separated);
\item every irreducible closed subset is the closure of a unique point;
\item the collection of quasi-compact opens forms a basis for the topology.
\end{enumerate}
Note that a theorem of Hochster (\cite{Hoc}) states that spectral spaces are exactly the topological spaces that are homeomorphic to (Zariski) spectra of rings, see \cite[Theorem 2.2]{Sch} for a more general statement. For any quasi-compact quasi-separated scheme $X$, the associated real scheme $X_r$ is a spectral space.

The class of \textbf{constructible subsets} of a spectral space is the smallest class of subsets stable by finite unions, finite intersections and complements which contains all quasi-compact open subsets.

\begin{definition}
\label{def:conscom}
If $M$ is a spectral space and $\Lambda$ is a ring, we define $D^b_{ctf}(M,\Lambda)$ to be the full subcategory of $D(M,\Lambda)$ consisting of complexes $A$ such that there is a finite increasing chain of closed constructible subsets of $M$
\begin{align}
\emptyset=Z_0\subset Z_1\subset\cdots\subset Z_n=M
\end{align}
such that for $i=1,\cdots,n$, $A_{|Z_i-Z_{i-1}}$ is the constant sheaf associated to a perfect complex of $\Lambda$-modules. In other words, these are complexes such that there exists a finite \emph{stratification} of $M$ into locally closed constructible subsets $M_i=Z_i-Z_{i-1}$ such that the restriction to each stratum is the constant sheaf associated to a perfect complex.
\end{definition}

\subsection{}
Similar to Definition~\ref{def:conscom}, recall that if $M$ is a spectral space and $\Lambda$ is a noetherian ring, a sheaf of $\Lambda$-modules $\mathcal{F}$ over $M$ is \textbf{constructible} if there is a finite stratification of $M$ into locally closed constructible subsets $M_i$ such that $\mathcal{F}_{|M_i}$ is the constant sheaf associated to a finitely generated $\Lambda$-module (\cite[Definition A.3]{Sch3}).
Then in the case where $\Lambda$ is noetherian, we can replace Definition~\ref{def:conscom} by the following equivalent characterizations:
\begin{proposition}
\label{prop:ctfequi}
Let $M$ be a spectral space and let $\Lambda$ be a noetherian ring. Then for $A\in D(M,\Lambda)$, the following are equivalent:
\begin{enumerate}
\item \label{num:ctf} $A\in D^b_{ctf}(M,\Lambda)$;
\item \label{num:boundedflat} $A$ can be represented by a bounded complex of constructible flat sheaves of $\Lambda$-modules;
\item \label{num:fintor} $A$ has finite $\operatorname{Tor}$-dimension and each $\mathcal{H}^i(A)$ is constructible.\footnote{Recall that a complex $A$ has $\operatorname{Tor}$-amplitude in $[a,b]$ if $H^n(A\overset{L}{\otimes}N)=0$ for any $n\notin[a,b]$ and any $\Lambda$-module $N$. A complex $A$ has finite $\operatorname{Tor}$-dimension if there exist $(a,b)$ such that $A$ has $\operatorname{Tor}$-amplitude in $[a,b]$.}
\end{enumerate}
\end{proposition}

\proof
It is clear that~\eqref{num:boundedflat} implies~\eqref{num:ctf}. To show that~\eqref{num:ctf} implies~\eqref{num:fintor}, by localization~\eqref{eq:locseq} we are reduced to the case where $A=Rj_!C$, where $j$ is the inclusion of a locally closed constructible subset and $C$ is a perfect complex of $\Lambda$-modules, in which case~\eqref{num:fintor} is clearly satisfied.

The equivalence between~\eqref{num:boundedflat} and~\eqref{num:fintor} is classical, whose proof can be easily adapted from the \'etale case, see \cite[Rapport 4.6]{SGA4.5}, \cite[Proposition 6.4.6]{Fu} and \cite[Tag 03TT]{Stack}.
\endproof

\subsection{}
The following result is an analogue of \cite[Theorem 6.3.10]{CD} for spectral spaces:
\begin{proposition}
\label{prop:conscomp}
Let $M$ be a spectral space of finite dimension and let $\Lambda$ be a ring. Then $D(M,\Lambda)$ is compactly generated, and the objects of $D^b_{ctf}(M,\Lambda)$ agree with the compact objects in $D(M,\Lambda)$.
\end{proposition}
\proof
The family of objects $Rj_!\Lambda$ where $j:U\to M$ is the inclusion of a quasi-compact open subset form a generating family of $D(M,\Lambda)$. By \cite[Tag 0902]{Stack}, any such $U$ is itself a spectral space of finite dimension. By \cite[Corollary 4.6]{Sch2}, the cohomological dimension of $U$ is bounded by the dimension of $M$. By \cite[Proposition 1.1.9]{CD}, the family $Rj_!\Lambda$ is a generating family of compact objects. It follows that the subcategory of compact objects in $D(M,\Lambda)$ agrees with the thick subcategory generated by those $Rj_!\Lambda$, and therefore lie in $D^b_{ctf}(M,\Lambda)$.

Conversely, we need to show that every object in $D^b_{ctf}(M,\Lambda)$ is compact. The localization sequence~\eqref{eq:locseq} implies that the functors $Rj_!$ and $j^*$ for $j$ the inclusion of a locally closed constructible subset preserve compact objects. The result then follows from the fact that perfect complexes of $\Lambda$-modules are compact objects in $D(\Lambda)$.
\endproof

\subsection{}
By \cite[Proposition 4.3.9]{BCR}, for any scheme $X$, the dimension of $X_r$ is bounded by the dimension of $X$. Therefore by Lemma~\ref{lm:mconscomp} and Proposition~\ref{prop:conscomp} we obtain the following result:
\begin{theorem}
\label{th:ctf=c}
Let $X$ be a quasi-compact quasi-separated scheme of finite dimension and let $\Lambda$ be a ring. Then the two subcategories $D^b_{ctf}(X_r,\Lambda)$ and $D_c(X_r,\Lambda)$ of $D(X_r,\Lambda)$ agree.
\end{theorem}

\begin{corollary}
\label{cor:prescons}
Let $\Lambda$ be a ring. Given a morphism $f:X\to S$ between quasi-compact quasi-separated schemes , the subcategory $D^b_{ctf}(X_r,\Lambda)$ of $D(X_r,\Lambda)$ is preserved by the following operations:
\begin{itemize}
\item $f^*$, $\overset{L}{\otimes}_S$, and $Rf_!$ for $f$ separated of finite type;
\item $Rf_*$ and $R\underline{Hom}_S$ for $f$ of finite type and $X$, $S$ quasi-excellent of finite dimension;
\item $f^!$ for $f$ separated of finite type and $X$, $S$ quasi-excellent of finite dimension.
\end{itemize}
\end{corollary}
\proof
By Theorem~\ref{th:ctf=c} we are reduced to show that these functors preserve $D_c(X_r,\Lambda)$. By~\ref{num:char0} we may assume that all schemes have characteristic $0$. In this case, by virtue of resolution of singularities, it is known that the six operations in the given situations preserve constructible objects in any motivic $\infty$-category, see \cite[Theorem 2.4.9]{BD} and \cite[Proposition 3.3]{DFKJ}.
\endproof

\begin{remark}
\label{rm:schfin}
In \cite[Remark 17.7.1]{Sch3}, Scheiderer raised the question whether the (higher) direct image functors of a morphism of finite type between excellent schemes preserve constructible sheaves on real schemes. Corollary~\ref{cor:prescons} gives a result in this direction by showing that the functor $Rf_*$ preserves $D^b_{ctf}$ for $f$ a morphism of finite type between quasi-excellent schemes of finite dimension.
\end{remark}

\subsection{}
\label{num:consfct}
In the rest of this section we deal with Grothendieck groups. For any stable $\infty$-category $\mathcal{C}$, its \textbf{Grothendieck group} $K_0(\mathcal{C})$ is the quotient of the free abelian group generated by its objects by the relations $[B]=[A]+[C]$ for all cofiber sequences $A\to B\to C\to A[1]$. For any ring $\Lambda$, let $K_0(\Lambda)$ be the Grothendieck group of the $\infty$-category of perfect complexes of $\Lambda$.\footnote{It is well-known that $K_0(\Lambda)$ agrees with the Grothendieck group of finitely generated projective $\Lambda$-modules (\cite[3.10]{TT}).}

If $M$ is a spectral space and $B$ is a ring, a $B$-valued \textbf{constructible function} is a function $\phi:M\to B$ such that there is a finite stratification of $M$ into locally closed constructible subsets $M_i$ such that $\phi_{|M_i}$ is constant for each $i$. This is equivalent to say that for each $m\in B$, $\phi^{-1}(\{m\})$ is locally closed and constructible, and is non-empty for only a finite number of $m$'s. We denote by $Cons(M, B)$ the $B$-algebra of constructible functions on $M$. As a $B$-module, $Cons(M,B)$ agrees with the free $B$-module generated by the characteristic functions of closed constructible subsets of $M$, and we have $Cons(M, B)\simeq Cons(M,\mathbb{Z})\otimes_{\mathbb{Z}}B$.

Let $M$ be a spectral space and let $\Lambda$ be a ring. For $A\in D^b_{ctf}(M,\Lambda)$ and $x\in M$, define the \textbf{local Euler-Poincar\'e index} $\chi(A)(x)\in K_0(\Lambda)$ as the class of the stalk $[A_x]\in K_0(\Lambda)$.\footnote{Note that when $\Lambda$ is a field, the class $[A_x]\in K_0(\Lambda)\simeq\mathbb{Z}$ agrees with $\sum_i(-1)^i\operatorname{dim}H^i(A_x)$.} We then have a canonical map
\begin{align}
\label{eq:locEP}
\begin{split}
\mathrm{Obj}(D^b_{ctf}(M,\Lambda))&\xrightarrow{}Cons(M,K_0(\Lambda))\\
A&\mapsto (x\mapsto\chi(A)(x)).
\end{split}
\end{align}
The map~\eqref{eq:locEP} factors through the Grothendieck group and induces a map
\begin{align}
\label{eq:locEPK0}
\chi:K_0(D^b_{ctf}(M,\Lambda))\xrightarrow{}Cons(M,K_0(\Lambda))
\end{align}
(see also \cite[Note 87$_1$]{ReS}). 
The following result is a variant of \cite[Theorem 9.7.1]{KS}:
\begin{proposition}
\label{prop:KS9.7.1}
The map~\eqref{eq:locEPK0} is an isomorphism.
\end{proposition}
\proof
For $a\in K_0(\Lambda)$, let $A$ be a perfect complex of $\Lambda$-modules whose class is $a$. If $j:Z\to M$ is the inclusion of a closed constructible subset, then map~\eqref{eq:locEPK0} sends $[Rj_!A]$ to $a_Z$. This shows that the map~\eqref{eq:locEPK0} is surjective.

Now prove the injectivity. Note that every element of $K_0(D^b_{ctf}(M,\Lambda))$ can be represented by a single complex in $D^b_{ctf}(M,\Lambda)$. Consider a complex $A\in D^b_{ctf}(M,\Lambda)$. Choose a finite stratification of $M$ into locally closed constructible subsets $M_i$ such that $A_{|M_i}$ is the constant sheaf associated to a perfect complex $C_i$ of $\Lambda$-modules. The localization sequence~\eqref{eq:locseq} shows that we have
\begin{align}
[A]=\sum_i[Rj_{i!}C_i]\in K_0(D^b_{ctf}(M,\Lambda))
\end{align}
where $j_i:M_i\to M$ is the inclusion. But to say that $\chi(A)=0$ means that for each $i$, $[C_i]=0\in K_0(\Lambda)$. This shows that $[A]=0$.
\endproof
From Proposition~\ref{prop:KS9.7.1} and Theorem~\ref{th:ctf=c} we deduce the following result:
\begin{corollary}
\label{cor:K0cons}
Let $X$ be a quasi-compact quasi-separated scheme of finite dimension and let $\Lambda$ be a ring. Then there is a canonical isomorphism
\begin{align}
K_0(D_c(X_r,\Lambda))\simeq Cons(X_r,K_0(\Lambda)).
\end{align}
\end{corollary}

\begin{remark}
\label{rm:SHret}
\begin{enumerate}
\item Similar arguments can be applied to the homotopy category of sheaves of spectra $\mathbf{SH}(X_r)$: by virtue of \cite[Corollary 3]{Bac}, the analogue of Theorem~\ref{th:ctf=c} states that the constructible objects of $\mathbf{SH}(X_r)$ are exactly the sheaves of spectra $C$ such that there exists a finite stratification of $X_r$ into locally closed constructible subsets such that the restriction of $C$ to each stratum is the constant sheaf associated to a compact spectrum. The analogue of Corollary~\ref{cor:K0cons} provides a canonical isomorphism
\begin{align}
K_0(\mathbf{SH}_c(X_r))\simeq Cons(X_r,\mathbb{Z}),
\end{align}
where we use the fact that the Grothendieck group of compact spectra agrees with that of finitely generated abelian groups, which follows from the canonical $t$-structure on spectra (\cite[Lemma 2]{Bac}).
\item Corollary~\ref{cor:K0cons} implies that any additive invariant on $D_c(X_r,\Lambda)$, such as the Euler characteristic or the characteristic class (\cite[Definition 5.1.3]{JY}), only depends on the constructible function defined by the local Euler-Poincar\'e index. This reflects the fact that the real \'etale site only lives on the characteristic $0$ fiber as in ~\ref{num:char0}, see \cite{Ill}.
\end{enumerate}
\end{remark}

\subsection{}
\label{num:CHM}
We now give an application in determining the Grothendieck group of the constructible rational motivic spectra. If $X$ is an excellent separated scheme of finite dimension, then by \cite[Proposition 2.10]{Bon2}, there is a bounded weight structure on $\mathbf{DM}_c(X,\mathbb{Q})$, called the \textbf{Chow weight structure}. We define the category of \textbf{Chow motives} with rational coefficients over $X$, denoted as $\mathbf{CHM}(X,\mathbb{Q})$, as the heart of the Chow weight structure.

\begin{remark}
As Bondarko's construction of the Chow weight structure is obtained by an abstract gluing procedure, the category $\mathbf{CHM}(X,\mathbb{Q})$ is quite mysterious in general. In some cases we have more explicit descriptions of this category:
\begin{enumerate}
\item By \cite[Corollary 2.12]{Bon2}, the category $\mathbf{CHM}(X,\mathbb{Q})$ contains the idempotent completion of the additive subcategory generated by elements of the form $p_*\mathbbold{1}_Y(d)[2d]$, where $d\in\mathbb{Z}$ and $p:Y\to X$ is a proper morphism with $Y$ regular. If $X$ is a scheme of finite type over an excellent separated scheme of dimension at most $2$, then the two categories coincide (\cite[Th\'eor\`eme 3.3]{Heb} and \cite[Theorem 2.1]{Bon2}).
\item By \cite[Theorem 3.17]{Jin}, if $X$ is quasi-projective over a perfect field, then $\mathbf{CHM}(X,\mathbb{Q})$ agrees with the category of Chow motives over $X$ defined in \cite[Definition 2.8]{CH}.
\end{enumerate}
\end{remark}

It is a general fact that weight structures have strong consequences on Grothendieck groups related to the heart.
For any additive category $\mathcal{A}$, denote by $K_0^{\oplus}(\mathcal{A})$ the quotient of the free abelian group generated by its objects by the relations $[B]=[A]+[C]$ if $B\simeq A\oplus C$. By \cite[Theorem 5.3.1]{Bon}, the inclusion $\mathbf{CHM}(X,\mathbb{Q})\to\mathbf{DM}_c(X,\mathbb{Q})$ induces an isomorphism
\begin{align}
\label{eq:K0DM}
K_0(\mathbf{DM}_c(X,\mathbb{Q}))\simeq K_0^{\oplus}(\mathbf{CHM}(X,\mathbb{Q})).
\end{align}

Combining~\eqref{eq:K0DM}, Corollary~\ref{cor:K0cons} and decomposition~\ref{eq:SHdec}, we deduce the following result:
\begin{corollary}
\label{cor:K0SHQ}
Let $X$ be an excellent separated scheme of finite dimension. Then there is a canonical isomorphism of abelian groups
\begin{align}
K_0(\mathbf{SH}_c(X,\mathbb{Q}))\simeq K_0^{\oplus}(\mathbf{CHM}(X,\mathbb{Q}))\oplus Cons(X_r,\mathbb{Z}).
\end{align}
\end{corollary}

\section{Generic base change}
\label{sec:genbc}
The main goal of this section is to prove the generic base change property for constructible complexes of real \'etale sheaves. The style of the proof is quite close to \cite[Th. finitude \S 2]{SGA4.5} and \cite[\S 2.4]{Cis}. See also \cite[Theorem 9.3.1]{Fu} for a detailed exposition. Throughout this section, we fix a ring $\Lambda$.

\subsection{}
\label{num:gbc}
Let $\mathcal{T}$ be a motivic $\infty$-category, or more generally any $\infty$-category fibered over schemes such that for any morphism of schemes $f:X\to Y$, the functor $f^*:\mathcal{T}(Y)\to\mathcal{T}(X)$ has a right adjoint $Rf_*:\mathcal{T}(X)\to\mathcal{T}(Y)$. Let $S$ be a scheme and let $f:X\to Y$ be a $S$-morphism of schemes. For an open subscheme $U$ over $S$ and an object $K\in\mathcal{T}(X)$, we say that $K$ \textbf{satisfies base change along $f$ over $U$} if the formation of $Rf_*K$ is compatible with any base change over $S$ which factors through $U$. In other words, for any morphism $p:V\to U$ with the Cartesian square
\begin{align}
\begin{split}
  \xymatrix@=10pt{
    X\times_SV \ar[r]^-{p_X} \ar[d]_-{f_V} & X\times_SU \ar[d]^-{f_U}\\
    Y\times_SV \ar[r]^-{p_Y} & Y\times_SU
  }
\end{split}
\end{align}
the canonical map $p_Y^*Rf_{U*}K_{|X\times_SU}\to Rf_{V*}p_X^*K_{|X\times_SU}$ is an isomorphism.

We say that an object $K\in\mathcal{T}(X)$ \textbf{satisfies generic base change along $f$ relatively to $S$} if there is an open dense subscheme $U$ of $S$ such that $K$ satisfies base change along $f$ over $U$.

\subsection{}
We start with a special case of the generic base change. If $\mathcal{C}$ is a closed symmetric monoidal $\infty$-category and $M\in\mathcal{C}$, we denote $M^\vee=R\underline{Hom}(M,\mathbbold{1})$ where $\mathbbold{1}$ is the monoidal unit. An object $M\in\mathcal{C}$ is called \textbf{dualizable} if the canonical map $M\otimes M^\vee\to R\underline{Hom}(M,M)$ is an isomorphism.
\begin{lemma}
\label{lm:cis2.4.4}
Let $\mathcal{T}$ be a motivic $\infty$-category, let $f:X\to S$ be a smooth morphism and let $K\in\mathcal{T}(X)$. Assume that $K$ is dualizable in $\mathcal{T}(X)$ and $Rf_!(K^\vee)$ is dualizable in $\mathcal{T}(S)$. Then for any Cartesian square of the form
\begin{align}
\begin{split}
  \xymatrix@=10pt{
    W \ar[r]^-{q} \ar[d]_-{g} & X \ar[d]^-{f}\\
    T \ar[r]^-{p} & S
  }
\end{split}
\end{align}
the canonical map $p^*Rf_*K\to Rg_*q^*K$ is an isomorphism.
\end{lemma}
\proof
The proof of \cite[Proposition 2.4.4]{Cis} works for any motivic $\infty$-category, using proper base change and relative purity. See Remark 2.4.5 of loc. cit.
\endproof

The following lemma slightly generalizes \cite[Lemma 2.4.9]{Cis}:
\begin{lemma}
\label{lm:cis2.4.9}
Let $S$ be a scheme and let $X\xrightarrow{f}Y\xrightarrow{j}Z$ be two composable $S$-morphisms, with $j$ and open immersion. Let $K\in\mathcal{T}(X)$ and let $U$ be an open subscheme of $S$. If $K$ satisfies base change along $jf$ over $U$, then $K$ satisfies base change along $f$ over $U$.
\end{lemma}
\proof
For a morphism $p:V\to U$, consider the Cartesian diagram
\begin{align}
\begin{split}
  \xymatrix@=10pt{
    X_V \ar[r]^-{p_X} \ar[d]_-{f_V} & X_U \ar[d]^-{f_U}\\
    Y_V \ar[r]^-{p_Y} \ar[d]_-{j_V} & Y_U \ar[d]^-{j_U} \\
    Z_V \ar[r]^-{p_Z} & Z_U.
  }
\end{split}
\end{align}
If $K$ satisfies base change along $jf$ over $U$, then for any such $p$ we have 
\begin{align}
p_Z^*Rj_{U*}Rf_{U*}K_{|X_U}\simeq Rj_{V*}Rf_{V*}p_X^*K_{|X_U}
\end{align}
Using the canonical identifications $j_U^*Rj_{U*}=id$ and $j_V^*Rj_{V*}=id$, we have
\begin{align}
\begin{split}
p_Y^*Rf_{U*}K_{|X_U}
&=
p_Y^*j_U^*Rj_{U*}Rf_{U*}K_{|X_U}
=
j_V^*p_Z^*Rj_{U*}Rf_{U*}K_{|X_U}\\
&\simeq
j_V^*Rj_{V*}Rf_{V*}p_X^*K_{|X_U}
=
Rf_{V*}p_X^*K_{|X_U}
\end{split}
\end{align}
and the result follows.
\endproof

\subsection{}
If a motivic $\infty$-category $\mathcal{T}$ satisfies some conditions on resolution of singularities (see \cite[2.4.1]{BD}, \cite[2.1.12]{JY}), then for any field $k$, every object in $\mathcal{T}_c(k)$ is dualizable (\cite[Remark 2.1.16]{JY}). Such conditions are automatically satisfied for $\mathcal{T}(X)=D(X_r,\Lambda)$ since by~\ref{num:char0} we may assume that all schemes have characteristic $0$. By the continuity property of $D(X_r,\Lambda)$ (\cite[Appendix A]{DFKJ}), we deduce the following result:
\begin{lemma}
\label{lm:gendual}
For every noetherian scheme $X$ and every constructible object $K$ of $D(X_r,\Lambda)$, there is a dense open subscheme $U$ of $X$ such that $K_{|U}$ is dualizable in $D(U_r,\Lambda)$.
\end{lemma}
From now on we focus on the case $\mathcal{T}(X)=D(X_r,\Lambda)$.
\begin{lemma}
\label{lm:fincons}
Let $f:X\to S$ be a finite morphism of schemes. 
Then the functor $Rf_*:D(X_r,\Lambda)\to D(S_r,\Lambda)$ is conservative.
\end{lemma}
\proof
Since $f$ is a finite morphism, then for every point $s\in S$ the fiber $f^{-1}(s)$ is finite, and the residue field of any point in the fiber is a finite extension of the residue field of $s$. It follows that the map $f_r:X_r\to S_r$ has finite discrete fibers, because for a finite field extension $k\to k'$ and an ordering $\sigma$ of $k$, there are only finitely many orderings on $k'$ extending the ordering $\sigma$ on $k$ (\cite[VIII 2.20]{Lam}). Therefore the functor $Rf_*:D(X_r,\Lambda)\to D(S_r,\Lambda)$ is clearly conservative.
\endproof

\subsection{}
Denote by $P(n)$ the following statement: for any integral quasi-excellent noetherian scheme of finite dimension $S$, and any open immersion with dense image $f:X\to Y$ between $S$-schemes of finite type such that the generic fiber of $X$ over $S$ has dimension at most $n$, any object $K\in D_c(X_r,\Lambda)$ satisfies generic base change along $f$ relatively to $S$.

\begin{proposition}
\label{prop:cis2.4.10}
Let $n\geqslant0$ be an integer, and assume that $P(n-1)$ holds. Let $S$ be an integral quasi-excellent noetherian scheme of finite dimension and let $f:X\to Y$ be a morphism between $S$-schemes of finite type such that $X$ is smooth over $S$ and the generic fiber of $X$ over $S$ has dimension $n$. Then any $K\in D(X_r,\Lambda)$ which is dualizable satisfies generic base change along $f$ relatively to $S$. 
\end{proposition}
\proof
Let $K\in D(X_r,\Lambda)$ be dualizable. The problem is local on $Y$, and therefore we may assume that $Y$ is affine over $S$. Choose a closed embedding $Y\to \mathbb{A}^d_S$ defined by $d$ functions $(g_k:Y\to\mathbb{A}^1_S)_{1\leqslant k\leqslant d}$. Then the generic fiber of $g_i\circ f$ has dimension at most $n-1$ (see \cite[7.5.3]{Fu}). For each $1\leqslant k\leqslant d$, applying the statement $P(n-1)$ to the morphism $f$ as a morphism of schemes over $\mathbb{A}^1_S$ via the structure morphism $g_k$, we know that there exists a dense open subscheme $U_k\subset\mathbb{A}^1_S$ such that $K$ satisfies base change along $f$ over $U_k$ via the structure morphism $g_k$. Let $V=\bigcup_{1\leqslant k\leqslant d}(g_k)^{-1}(U_k)$ and denote by $j:V\to Y$ the canonical open embedding. Then $K_{|f^{-1}(V)}$ satisfies base change along $f_V:f^{-1}(V)\to V$ over $S$. Let $T=Y-V$ be the complement of $V$ (with any scheme structure). Then $T$ satisfies
\begin{align}
T\subset Y\cap\left((\mathbb{A}^1_S-U_1)\times_S\cdots\times_S (\mathbb{A}^1_S-U_d)\right).
\end{align}
Since the generic fiber of $(\mathbb{A}^1_S-U_1)\times_S\cdots\times_S (\mathbb{A}^1_S-U_d)\to S$ is finite, by shrinking $S$, we may assume that $T$ is finite over $S$.

Let $\bar{Y}$ be the closure of $Y$ in $\mathbb{P}^n_S$. 
The image $S_1$ of $\bar{Y}-Y\to S$ is closed in $S$ (since $\bar{Y}$ is proper over $S$) and does not contain the generic point (since the generic fiber of $X$ is non-empty), therefore by replacing $S$ by $S-S_1$, we may assume that $\bar{Y}-V$ is also finite over $S$. 
By Lemma~\ref{lm:cis2.4.9}, we may replace $Y$ by $\bar{Y}$, which amounts to say that we may assume that the structure morphism $p:Y\to S$ is proper. Therefore we are reduced to the following situation: $Y$ is proper over $S$, and there exists a open immersion $j:V\to Y$ with dense image whose complement $Y-V$ finite over $S$ such that $K_{|f^{-1}(V)}$ satisfies base change along $f_V:f^{-1}(V)\to V$ over $S$, that is, the formation of $Rj_!j^*Rf_*K$ is compatible with any base change over $S$. 

Denote by $i:Y-V\to Y$ the canonical closed immersion. We have the localization sequence
\begin{align}
\label{eq:locflower*}
Rj_!j^*Rf_*K\to Rf_*K\to Ri_*i^*Rf_*K\to Rj_!j^*Rf_*K[1].
\end{align}
Therefore it remains to show that after shrinking $S$, the formation of $Ri_*i^*Rf_*K$ is compatible with any base change over $S$. Since $Ri_*=Ri_!$ commutes with any base change, it suffices to prove that after shrinking $S$, the formation of $i^*Rf_*K$ is compatible with any base change over $S$. Since the composition $pi$ is finite, by Lemma~\ref{lm:fincons}, it suffices to prove that after shrinking $S$, the formation of $Rp_*Ri_*i^*Rf_*K$ is compatible with any base change over $S$. Applying the functor $Rp_*$ to~\eqref{eq:locflower*}, we obtain the following distinguished triangle:
\begin{align}
\label{eq:locflower*p}
Rp_*Rj_!j^*Rf_*K\to Rp_*Rf_*K\to Rp_*Ri_*i^*Rf_*K\to Rp_*Rj_!j^*Rf_*K[1].
\end{align}
Since the composition $pf:X\to S$ is smooth over $S$, we may apply Lemma~\ref{lm:cis2.4.4} and Lemma~\ref{lm:gendual} to obtain that, after shrinking $S$, the formation of $Rp_*Rf_*K=R(pf)_*K$ is compatible with any base change over $S$. On the other hand, since $p$ is proper and the formation of $Rj_!j^*Rf_*K$ is compatible with any base change over $S$, we deduce from the proper base change that the formation of $Rp_*Rj_!j^*Rf_*K$ is compatible with any base change over $S$. We conclude using the distinguished triangle~\eqref{eq:locflower*p}.
\endproof

\begin{theorem}
\label{th:ret_genbc}
Let $S$ be a quasi-excellent noetherian scheme of finite dimension and let $f:X\to Y$ be a morphism between $S$-schemes of finite type. Then every object of $D_{c}(X_r,\Lambda)$ (respectively $\mathbf{SH}_c(X_r)$) satisfies generic base change along $f$ relatively to $S$.
\end{theorem}
\proof
We prove the case of $D_{c}(X_r,\Lambda)$, the case of $\mathbf{SH}_c(X_r)$ being very similar. The problem is local on $Y$, so we may assume that $Y$ is affine. By virtue of the hypercohomology spectral sequence for a finite open affine cover of $X$ (which holds since the real \'etale $\infty$-topos hypercomplete, see \cite[Theorem B.13]{ES}), we reduce the problem to the case where $X$ is also affine. So we may assume that there exists a compactification $f=pj$ with $p$ proper and $j$ open immersion. By proper base change, we may assume that $f$ is an open immersion with dense image. By working with each irreducible component of $S$, we may also assume that $S$ is integral. Therefore we are reduced to prove $P(n)$.

We use induction on $n$. The case $n=-1$ is clear, and we may assume that $n\geqslant0$ and $P(n-1)$ holds. By~\ref{num:char0}, we may assume that all schemes have characteristic $0$. Then there exists a dense open subscheme $S_0$ of $S$ and a dense open subscheme $U$ of $X$ which is smooth over $S_0$. Replacing $S$ by $S_0$, we may assume that there is an open immersion $j:U\to X$ with dense image such that $U$ is smooth over $S$.

Let $K\in D_{c}(X_r,\Lambda)$. By Lemma~\ref{lm:gendual}, by shrinking $U$, we may also assume that $K_{|U}$ is dualizable in $D(U_r,\Lambda)$. Let $i:Z\to X$ be the closed complement (with any scheme structure). 
Consider the following localization sequence:
\begin{align}
\label{eq:locflower*p}
R(fi)_*i^!K\to Rf_*K\to R(fj)_*j^*K\to R(fi)_*i^!K[1].
\end{align}
By Corollary~\ref{cor:prescons}, the object $i^!K$ lies in $D_{c}(Z_r,\Lambda)$.
Applying Proposition~\ref{prop:cis2.4.10} to $K$ along the morphism $fj$ and the statement $P(n-1)$ to $i^!K$ along the morphism $fi$, we know that there exists a dense open subscheme $W$ of $S$ such that the formations of $R(fi)_*i^!K$ and $R(fj)_*j^*K$ are compatible with any base change over $S$ which factors through $W$. Therefore the same property holds for $Rf_*K$, which finishes the proof.
\endproof

\begin{corollary}
\label{cor:genbcsh}
Let $S$ be a quasi-excellent noetherian scheme of finite dimension and let $f:X\to Y$ be a morphism between $S$-schemes of finite type. Then every object of $\mathbf{SH}_c(X,\mathbb{Q})$ satisfies generic base change along $f$ relatively to $S$.
\end{corollary}
\proof
By the decomposition~\eqref{eq:SHdec}, we only need to prove the generic base change for $D_c(X_r,\mathbb{Q})$ and $\mathbf{DM}_c(X,\mathbb{Q})$. The first case is Theorem~\ref{th:ret_genbc}. For the second case, by \cite[Theorem 5.2.2]{CD} we know that $\mathbf{DM}_c(X,\mathbb{Q})$ agrees with constructible $h$-motives with rational coefficients, for which generic base change is proved in \cite[Theorem 2.4.2]{Cis}.
\endproof

\subsection{}
Recall that the $b$-topology is obtained by gluing the \'etale topology and the real \'etale topology (\cite[Definition 2.3]{Sch3}). Denote by $j$ (respectively $i$) the canonical inclusion of the \'etale topos (respectively the real \'etale topos) into the $b$ topos. Then we have the following generic base change result for complexes of $b$-sheaves:
\begin{corollary}
\label{cor:genbcb}
Let $S$ be a quasi-excellent noetherian scheme of finite dimension and let $f:X\to Y$ be a morphism between $S$-schemes of finite type. Let $\Lambda$ be a noetherian ring, and let $K\in D(X_b,\Lambda)$ be such that 
\begin{enumerate}
\item $i^*K\in D_{c}(X_r,\Lambda)$;
\item $j^*K\in D^b_{ctf}(X_{\textrm{\'et}},\Lambda)$;
\item There exists an integer $n\in\mathcal{O}^*(X)$ such that $n\cdot j^*K=0$.
\end{enumerate}
Then $K$ satisfies generic base change along $f$ relatively to $S$.
\end{corollary}
\proof
By virtue of \cite[Remarks 16.1]{Sch3}, the result follows from the real \'etale case and the \'etale case, which follow from Theorem~\ref{th:ret_genbc} and \cite[Theorem 2.4.2]{Cis} combined with \cite[Theorem 6.3.11]{CD} respectively.
\endproof

\end{document}